\documentclass[final,1p,times]{elsarticle}

\usepackage{amssymb}
 \usepackage{amsthm}
\usepackage{amscd}
\usepackage{amsmath}
\usepackage{amsfonts}
\usepackage{amssymb}
\usepackage{graphicx}
\usepackage{amsfonts}
\usepackage{indentfirst,latexsym,bm,amssymb,amsthm,amsmath,txfonts}
\newtheorem{theorem}{Theorem}[section]

\newtheorem{Lem}[theorem]{Lemma}

\newtheorem{Thm}[theorem]{Theorem}
\numberwithin{equation}{section}

\newcommand{\ra}{\rightarrow}

\newcommand{\f}{\frac}

\newcommand{\be}{\begin{equation}}
\renewcommand{\ra}{\rightarrow}
\newcommand{\ee}{\end{equation}}
\newcommand{\bea}{\begin{eqnarray}}
\newcommand{\eea}{\end{eqnarray}}
\newcommand{\bna}{\begin{eqnarray*}}
\newcommand{\ena}{\end{eqnarray*}}

\renewcommand{\le}{\left}
\newcommand{\ri}{\right}

\journal{***}

\begin{document}

\begin{frontmatter}

\title{Trudinger-Moser embedding on the hyperbolic space}

 \author{Yunyan Yang}
 \ead{yunyanyang@ruc.edu.cn}
 \author{Xiaobao Zhu}
 \ead{zhuxiaobao@ruc.edu.cn}
\address{ Department of Mathematics,
Renmin University of China, Beijing 100872, P. R. China}

\begin{abstract}

Let $(\mathbb{H}^n,g)$ be the hyperbolic space of dimension $n$. By
our previous work (Theorem 2.3 of \cite{YangJFA}), for any
$0<\alpha<\alpha_n$, there exists a constant $\tau>0$ depending only
on $n$ and $\alpha$ such that
 \be\label{0.1}\sup_{u\in
W^{1,n}(\mathbb{H}^n),\,\|u\|_{1,\tau}\leq 1}
 \int_{\mathbb{H}^n}\le(e^{\alpha|u|^{\f{n}{n-1}}}-\sum_{k=0}^{n-2}\f{\alpha^k|u|^{\f{nk}{n-1}}}{k!}\ri)dv_g
 <\infty,\ee
 where $\alpha_n=n\omega_{n-1}^{\f{1}{n-1}}$, $\omega_{n-1}$ is the area of the unit sphere $\mathbb{S}^n$,
 and
 $\|u\|_{1,\tau}=\|\nabla_gu\|_{L^n(\mathbb{H}^n)}+\tau\|u\|_{L^n(\mathbb{H}^n)}$.
 In this note we shall improve (\ref{0.1}). Particularly we show that for any $0<\alpha<\alpha_n$ and any $\tau>0$,
  (\ref{0.1}) holds with
 the definition of $\|u\|_{1,\tau}$ replaced by
 $\le(\int_{\mathbb{H}^n}(|\nabla_gu|^n+\tau|u|^n)dv_g\ri)^{1/n}$.
  We solve this problem by
 gluing local uniform estimates.
\end{abstract}


\end{frontmatter}

\section{Introduction}
 Let $\Omega$ be a bounded smooth domain in $\mathbb{R}^n$. The classical Trudinger-Moser inequality
\cite{Moser,Pohozaev,Trudinger} says \be\label{T-M-Omega}\sup_{u\in
W_0^{1,n}(\Omega),\,\|u\|_{W_0^{1,n}(\Omega)}\leq 1}\int_\Omega
e^{\alpha_n|u|^{\f{n}{n-1}}}dx\leq C|\Omega|\ee for some constant
$C$ depending only on $n$, where $W_0^{1,n}(\Omega)$ is the usual
Sobolev space and $|\Omega|$ denotes the Lebesgue measure of
$\Omega$. In the case $\Omega$ is an unbounded domain of
$\mathbb{R}^n$, the above integral is infinite, but it was shown by
  Cao \cite{Cao}, Panda \cite{Panda} and do
 \'O \cite{doo1} that for any $\tau>0$ and any $\alpha<\alpha_n$ there
 holds
\be\label{sub-eucl} \sup_{u\in
W^{1,n}(\mathbb{R}^n),\,\,\int_{\mathbb{R}^{n}}(|\nabla
    u|^{n}+\tau|u|^{n})dx\leq1}\int_{\mathbb{R}^{n}}\le(e^{\alpha|u|^{\frac{n}{n-1}}}-\sum_{k=0}^{n-2}
    \frac{\alpha^{k}|u|^{\frac{nk}{n-1}}}{k!}\ri)dx<\infty.\ee
  Later Ruf \cite{Ruf}, Li-Ruf \cite{LiRuf} and
 Adimurthi-Yang \cite{Adi-Yang} obtained (\ref{sub-eucl}) in the
 critical case $\alpha=\alpha_n$.

The study of Trudinger-Moser inequalities on compact Riemannian
manifolds can be traced back to Aubin \cite{Aubin10}, Cherrier
\cite{cherrier1,cherrier2}, and Fontana \cite{Fontana}. A particular
case is as follows. Let $(M,g)$ be an $n$-dimensional compact
Riemannian manifold without boundary. Then there holds
\be\label{T-M-manifold}\sup_{\int_M|\nabla_gu|^ndv_g\leq
1,\,\int_Mudv_g=0}\int_Me^{\alpha_n|u|^{\f{n}{n-1}}}dv_g<\infty.\ee

 In view of (\ref{sub-eucl}), it is natural to consider extension of (\ref{T-M-manifold}) on
 complete noncompact Riemannian manifolds. In \cite{YangJFA} we
 obtained the following results:  Let $(M,g)$ be a complete noncompact Riemannian
 manifold. If the Trudinger-Moser
 inequality holds on it, then there holds
 $\inf_{x\in M}{\rm vol}_g(B_1(x))>0$. If the Ricci curvature has lower bound, say ${\rm Ric}_g(M)\geq
 -K$, the injectivity radius has a positive lower bound $i_0$, then
 for any $\alpha<\alpha_n$ there exists a constant $\tau>0$
 depending only on $\alpha$, $n$, $K$, and $i_0$ such that
 \be\label{t-m-manifold}\sup_{\le(\int_M|\nabla u|^ndv_g\ri)^{1/n}+\tau\le(\int_M|u|^ndv_g\ri)^{1/n}\leq 1}
 \int_M\le(e^{\alpha|u|^{\frac{n}{n-1}}}-\sum_{k=0}^{n-2}
    \frac{\alpha^{k}|u|^{\frac{nk}{n-1}}}{k!}\ri)dv_g<\infty.\ee
   Since $\tau$ depends on $\alpha$, (\ref{t-m-manifold}) is weaker than (\ref{sub-eucl}) when $(M,g)$ is replaced by
   $\mathbb{R}^n$. Moreover, the condition that ${\rm Ric}_g(M)$ has
   lower bound is not necessary for the validity of the Trudinger-Moser
 inequality.

 In this note, we shall improve (\ref{t-m-manifold}) in a special case that $(M,g)$ is the hyperbolic
 space $(\mathbb{H}^n,g)$, a simply connected Riemannian
 manifold with constant sectional curvature $-1$. Particularly we
 have the following:

\begin{Thm}\label{main theorem} Let $(\mathbb{H}^n,g)$ be an
$n$-dimensional hyperbolic space,
$\alpha_{n}=n\omega_{n-1}^{\frac{1}{n-1}}$, where $\omega_{n-1}$ is
the measure of the unit sphere in $\mathbb{R}^{n}$. Then for any
$\alpha<\alpha_{n}$, any $\tau>0$, and any $u\in
W^{1,n}(\mathbb{H}^n)$ satisfying $\int_{\mathbb{H}^n}(|\nabla_g\,
u|^{n}+\tau|u|^{n})dv_{g}\leq 1$, there exists some constant $\beta$
depending only on $n$ and $\tau$ such that

\be\label{e1-1} \int_{\mathbb{H}^{n}}
    \left(e^{\alpha|u|^{\frac{n}{n-1}}}-\sum_{k=0}^{n-2}
    \frac{\alpha^{k}|u|^{\frac{nk}{n-1}}}{k!}\right)dv_{g}\leq \beta.\ee
\end{Thm}

The proof of Theorem \ref{main theorem} is based on local uniform
estimates  (Lemma 2.1 below). This idea comes from \cite{YangJFA}
and was also used in \cite{Y-Heisenberg,YZ-2012}. The remaining part
of this note is organized as follows. In Section 2 we derive local
uniform Trudinger-Moser inequalities; In Section 3, Theorem 1.1 is
proved.

\section{Local estimates}
To get (\ref{e1-1}), we need the following uniform local estimates
which is an analogy of (\cite{Y-Heisenberg}, Lemma 4.1) or
(\cite{YZ-2012}, Lemma 1), and of its own interest.

\begin{Lem}\label{Trudinger-Moser ineq. on bounded domain}
 For any $p\in\mathbb{H}^n$, any $R>0$, and any
 $u\in W_{0}^{1,n}({B}_{R}(p))$ with $\int_{{B}_{R}(p)}|\nabla_g u|^{n}dv_{g}\leq1$, there exists some
 constant $C_n$ depending only on $n$ such that
\begin{align}\label{e2-1}
  \int_{{B}_{R}(p)}\left(e^{\alpha_{n}|u|^{\frac{n}{n-1}}}
                -\sum_{k=0}^{n-2}\frac{\alpha_{n}^{k}|u|^{\frac{nk}{n-1}}}{k!}\right)dv_{g}
                \leq C_n(\sinh R)^{n}\int_{{B}_{R}(p)}|\nabla_g
                u|^{n}dv_{g},
\end{align}
where ${B}_R(p)$ denotes the geodesic ball of $(\mathbb{H}^n,g)$
which is centered at $p$ with radius $R$.
\end{Lem}

{\it Proof}.
 It is well known, see for example \cite{Chavel-book}, II.5, Theorem
 1, that there exists a homomorphism $\varphi:\mathbb{H}^n\ra
 D=\{x\in\mathbb{R}^{n}:|x|<1\}$ such that $\varphi(p)=0$,  that in these coordinates the
 Riemannian metric $g$ can be represented by
 $$g(x)=\frac{4}{(1-|x|^{2})^{2}}g_{0}(x),$$ where
$g_{0}(x)=\sum_{i=1}^{n}(dx^{i})^{2}$ is the standard Euclidean
metric on $\mathbb{R}^{n}$, and that
$$\varphi(B_R(p))=\mathbb{B}_{\tanh \f{R}{2}}(0),$$
where $\mathbb{B}_{r}(0)\subset\mathbb{R}^n$ denotes a ball centered
at $0$ with radius $r$. Moreover, the corresponding polar
coordinates $(r,\theta) \in[0,\infty)\times \mathbb{S}^{n-1}$ reads
$$g=dr^2+(\sinh r)^2d\theta^2,$$
where $d\theta^2$ is the standard metric on $\mathbb{S}^{n-1}$.

    Denote $f=\frac{2}{1-|x|^{2}}$, then $g=f^{2}g_{0}$, $|\nabla_g u|=f^{-1}|\nabla_{g_0}(u\circ\varphi^{-1})|$ and
$dv_{g}=f^{n}dv_{g_{0}}$. Calculating directly, we have
    \begin{align}\label{e2}
     \int_{{B}_{R}(0)}|\nabla_g u|^{n}dv_{g}
    =\int_{\mathbb{B}_{\tanh \f{R}{2}}(0)}|\nabla_{g_0}(u\circ\varphi^{-1})|^{n}dv_{g_{0}}.
    \end{align}
Since $u\in W_{0}^{1,n}({B}_{R}(p))$, we have $u\circ\varphi^{-1}\in
W_{0}^{1,n}(\mathbb{B}_{\tanh\f{R}{2}}(0))$. Noting that
$\int_{{B}_{R}(p)}|\nabla_g u|^{n}dv_{g}\leq1$, we have by
(\ref{e2})
$$\int_{\mathbb{B}_{\tanh \f{R}{2}}(0)}|\nabla_{g_0}(u\circ\varphi^{-1})|^{n}dv_{g_{0}}\leq 1.$$
The standard Trudinger-Moser inequality (\ref{T-M-Omega})
implies
 \bna
 \int_{\mathbb{B}_{\tanh\f{R}{2}}(0)}\left(e^{\alpha_{n}|u\circ\varphi^{-1}|^{\frac{n}{n-1}}}
                -\sum_{k=0}^{n-2}\frac{\alpha_{n}^{k}|u\circ\varphi^{-1}|^{\frac{nk}{n-1}}}{k!}\right)dv_{g_{0}}
                &=&\int_{\mathbb{B}_{\tanh\f{R}{2}}(0)}\sum_{k=n-1}^\infty\f{\alpha_n^k|u\circ\varphi^{-1}|^{\f{nk}{n-1}}}{k!}dv_{g_0}\\
                &\leq&\int_{\mathbb{B}_{\tanh\f{R}{2}}(0)}\sum_{k=n-1}^\infty\f{\alpha_n^k|\f{u\circ\varphi^{-1}}
                {\|\nabla_{g_0}(u\circ\varphi^{-1})\|_{L^n}}|^{\f{nk}{n-1}}}{k!}dv_{g_0}\\
                &&\quad\times\int_{\mathbb{B}_{\tanh\f{R}{2}}(0)}|\nabla_{g_0}
                (u\circ\varphi^{-1})|^{n}dv_{g_{0}}\\
                &\leq& C_n\le(\tanh\f{R}{2}\ri)^n\int_{\mathbb{B}_{\tanh\f{R}{2}}(0)}|\nabla_{g_0}
                (u\circ\varphi^{-1})|^{n}dv_{g_{0}},
\ena where $C_n$ is a constant depending only on $n$.
    This together with (\ref{e2}) immediately leads to
     \bea\nonumber\int_{{B}_{R}(p)}\left(e^{\alpha_{n}|u|^{\frac{n}{n-1}}}
                -\sum_{k=0}^{n-2}\frac{\alpha_{n}^{k}|u|^{\frac{nk}{n-1}}}{k!}\right)dv_{g}&=&
                \int_{\mathbb{B}_{\tanh\f{R}{2}}(0)}\left(e^{\alpha_{n}|u\circ\varphi^{-1}|^{\frac{n}{n-1}}}
                -\sum_{k=0}^{n-2}\frac{\alpha_{n}^{k}|u\circ\varphi^{-1}|^{\frac{nk}{n-1}}}{k!}\right)f^ndv_{g_{0}}\\
                \nonumber&\leq&C_n\left(\frac{2\tanh\f{R}{2}}{1-\le(\tanh\f{R}{2}\ri)^{2}}\right)^{n}\int_{\mathbb{B}_{\tanh\f{R}{2}}(0)}|\nabla_{g_0}
                (u\circ\varphi^{-1})|^{n}dv_{g_{0}}\\\label{abov}&=&C_n(\sinh R)^{n}\int_{{B}_{R}(p)}|\nabla_g
                u|^{n}dv_{g}.\eea
This is exactly (\ref{e2-1}) and thus ends the proof of the lemma.   $\hfill\Box$\\

As a corollary of Lemma 2.1, the following estimates can be compared
with (\ref{T-M-Omega}).\\

\noindent {\bf Corollary 2.2.} {\it For any $p\in\mathbb{H}^n$, any
$R>0$, and any $u\in W_0^{1,n}(B_R(p))$ with
$\int_{B_R(p)}|\nabla_gu|^ndv_g\leq 1$, there exists some constant
$C$ depending only on $n$ such that
 \be
 \label{hy}
 \f{1}{{\rm
 Vol}_g(B_R(p))}\int_{B_R(p)}e^{\alpha_n|u|^{\f{n}{n-1}}}dv_g\leq
 C\f{\sinh R}{R}.
\ee
}

 {\it Proof}. Since
 $$\lim_{R\ra 0+}\f{{\rm Vol}_g(B_R(p))}{R(\sinh R)^{n-1}}=\lim_{R\ra\infty}\f{{\rm Vol}_g(B_R(p))}{R(\sinh R)^{n-1}}=1,$$
 it follows from (\ref{abov}) that there exists some constant
 $C$ depending only on $n$ such that
 \be\label{sb}\f{1}{{\rm
 Vol}_g(B_R(p))}\int_{{B}_{R}(p)}\left(e^{\alpha_{n}|u|^{\frac{n}{n-1}}}
                -\sum_{k=0}^{n-2}\frac{\alpha_{n}^{k}|u|^{\frac{nk}{n-1}}}{k!}\right)dv_{g}
                \leq C\f{\sinh R}{R}.\ee
 In particular,
 $$\int_{{B}_{R}(p)}|u|^{n}dv_g\leq  C\f{\sinh R}{R}{\rm
 Vol}_g(B_R(p)).$$
 Here and in the sequel we often denote various constants by
 the same $C$, the reader can easily distinguish them from the context.
 Noting that for any $q$, $0\leq q\leq n$,
 $$\int_{{B}_{R}(p)}|u|^qdv_g\leq {\rm
 Vol}_g(B_R(p))+\int_{{B}_{R}(p)}|u|^ndv_g,$$
 we conclude
 \bea\label{res}\int_{{B}_{R}(p)}\sum_{k=0}^{n-2}\frac{\alpha_{n}^{k}|u|^{\frac{nk}{n-1}}}{k!}
 dv_g\leq C\f{\sinh R}{R}{\rm
 Vol}_g(B_R(p)).\eea
 Combining (\ref{sb}) and (\ref{res}), we obtain (\ref{hy}). $\hfill\Box$

\section{Proof of Theorem {\ref{main theorem}}}
 In this section, we will prove Theorem \ref{main theorem} by gluing local estimates
 (\ref{e2-1}).\\

{\it Proof of Theorem \ref{main theorem}}. Let $R$ be a positive
real number which will be determined later. By (\cite{Hebey}, Lemma
1.6) we
 can find a sequence of points
$\{x_{i}\}_{i=1}^{\infty}\subset\mathbb{H}^{n}$ such that
$\cup_{i=1}^{\infty}{B}_{\frac{R}{2}}(x_{i})=\mathbb{H}^{n}$, that
${B}_{\frac{R}{4}}(x_{i})\cap {B}_{\frac{R}{4}}(x_{j})=\varnothing$
for any $ i\neq j$, and that for any $x\in\mathbb{H}^{n}$, $x$
belongs to at most $N$ balls ${B}_{R}(x_{i})$, where $N$ depends
only on $n$. Let $\phi_{i}$ be the cut-off function satisfies the
following conditions: $(i)$ $\phi_{i}\in
C_{0}^{\infty}({B}_{R}(x_{i}))$;
   $(ii)$ $0\leq\phi_{i}\leq1$ on ${B}_R(x_i)$ and $\phi_{i}\equiv 1$ on ${B}_{{R}/{2}}(x_{i})$;
   $(iii)$ $|\nabla_g \phi_{i}(x)|\leq{4}/{R}$.
Let $\tau>0$ be fixed. For any $u\in W^{1,n}(\mathbb{H}^{n})$
satisfying
\begin{align}\label{e3-1}
\int_{\mathbb{H}^{n}}(|\nabla_g u|^{n}+\tau|u|^{n})dv_{g}\leq1,
\end{align}
we have $\phi_{i}u\in W_{0}^{1,n}({B}_{R}(x_{i}))$. For any
$\epsilon>0$, using an elementary inequality $ab\leq \epsilon
a^2+\f{1}{4\epsilon}b^2$, we find some constant $C$ depending only
on $n$ and $\epsilon$ such that
\begin{align}\label{e3-2}
\int_{{B}_{R}(x_{i})}|\nabla_g(\phi_{i}u)|^{n}dv_{g}\leq&(1+\epsilon)\int_{{B}_{R}(x_{i})}\phi_{i}^{n}
|\nabla_g u|^{n}dv_{g}
         +C\int_{{B}_{R}(x_{i})}|\nabla_g \phi_{i}|^{n}|u|^{n}dv_{g}\nonumber\\
         \leq&(1+\epsilon)\int_{{B}_{R}(x_{i})}|\nabla_g u|^{n}dv_{g}
         +\frac{4^nC}{R^{n}}\int_{{B}_{R}(x_{i})}|u|^{n}dv_{g}\nonumber\\
         \leq&(1+\epsilon)\int_{{B}_{R}(x_{i})}(|\nabla_g u|^{n}+\tau|u|^{n})dv_{g},
\end{align}
where in the last inequality we choose a sufficiently large $R$ to
make sure $\frac{4^nC}{R^{n}}\leq(1+\epsilon)\tau$. Let
$\alpha_{\epsilon}=\frac{\alpha_{n}}{(1+\epsilon)^{1/(n-1)}}$ and
$\widetilde{\phi_{i}u}=\frac{\phi_{i}u}{(1+\epsilon)^{1/n}}$. Noting
that $\widetilde{\phi_{i}u}\in W_{0}^{1,n}({B}_{R}(x_{i}))$, we have
by (\ref{e3-2}) and Lemma \ref{Trudinger-Moser ineq. on bounded
domain}
\begin{align}\label{e3-3}
    \int_{{B}_{\frac{R}{2}}(x_{i})}\left(e^{\alpha_{\epsilon}|u|^{\frac{n}{n-1}}}-\sum_{k=0}^{n-2}\frac{\alpha_{\epsilon}^{k}|u|^{\frac{nk}{n-1}}}{k!}\right)dv_{g}
\leq&\int_{{B}_{R}(x_{i})}\left(e^{\alpha_{\epsilon}|\phi_{i}u|^{\frac{n}{n-1}}}-\sum_{k=0}^{n-2}\frac{\alpha_{\epsilon}^{k}|\phi_{i}u|^{\frac{nk}{n-1}}}{k!}\right)dv_{g}\nonumber\\
   =&\int_{{B}_{R}(x_{i})}\left(e^{\alpha_{n}|\widetilde{\phi_{i}u}|^{\frac{n}{n-1}}}-\sum_{k=0}^{n-2}\frac{\alpha_{n}^{k}|\widetilde{\phi_{i}u}|^{\frac{nk}{n-1}}}{k!}\right)dv_{g}
   \nonumber\\
\leq&C_n(\sinh  R)^n\int_{{B}_{R}(x_{i})}|\nabla_g(\widetilde{\phi_{i}u})|^{n}dv_{g}\nonumber\\
\leq&C(\sinh R)^n\int_{{B}_{R}(x_{i})}(|\nabla_g
u|^{n}+\tau|u|^{n})dv_{g},
\end{align}
where $C$ is a constant depending only on $n$ and $\tau$. By the
choice of $\{x_i\}_{i=1}^\infty$ and (\ref{e3-3}), we have
\begin{align}\label{e3-4}
   \int_{\mathbb{H}^{n}}\left(e^{\alpha_{\epsilon}|u|^{\frac{n}{n-1}}}-\sum_{k=0}^{n-2}\frac{\alpha_{\epsilon}^{k}|u|^{\frac{nk}{n-1}}}{k!}\right)dv_{g}
\leq&\int_{\cup_{i=1}^{\infty}{B}_{\frac{R}{2}}(x_{i})}\left(e^{\alpha_{\epsilon}|u|^{\frac{n}{n-1}}}-\sum_{k=0}^{n-2}\frac{\alpha_{\epsilon}^{k}|u|^{\frac{nk}{n-1}}}{k!}\right)dv_{g}\nonumber\\
\leq&\sum_{i=1}^{\infty}\int_{{B}_{\frac{R}{2}}(x_{i})}\left(e^{\alpha_{\epsilon}|u|^{\frac{n}{n-1}}}-\sum_{k=0}^{n-2}\frac{\alpha_{\epsilon}^{k}|u|^{\frac{nk}{n-1}}}{k!}\right)dv_{g}\nonumber\\
\leq&\sum_{i=1}^{\infty}C(\sinh R)^{n}\int_{{B}_{R}(x_{i})}(|\nabla_g u|^{n}+\tau|u|^{n})dv_{g}\nonumber\\
\leq&CN(\sinh R)^{n}\int_{\mathbb{H}^{n}}(|\nabla_g u|^{n}+\tau|u|^{n})dv_{g}\nonumber\\
\leq&CN(\sinh R)^{n}
\end{align}
for some constant $C$ depending only on $n$ and $\tau$. For any
$\alpha<\alpha_n$, we can choose $\epsilon>0$ sufficiently small
such that $\alpha<\alpha_\epsilon$. This ends the proof of
Theorem \ref{main theorem}. $\hfill\Box$\\

{\bf Acknowledgement.} This work is supported by the NSFC
 11171347.

\end{document}